\documentclass[12pt]{elsarticle}

\usepackage{amsfonts,amsmath,amsthm,amscd,amssymb,latexsym}
\usepackage[english]{babel}

\theoremstyle{plain}

\newtheorem{thm}{Theorem}
\newtheorem{lem}[thm]{Lemma}
\newtheorem{prop}[thm]{Proposition}
\newtheorem{cor}[thm]{Corollary}
\newdefinition{rmk}{Remark}
\newdefinition{dfn}{Definition}
\newdefinition{exm}{Example}
\newdefinition{question}{Question}
\newproof{pf}{Proof}

\newtheorem{thmm}{Theorem}

\begin{document}

\title{Functionally $\sigma$-discrete mappings\\ and a generalization of Banach's theorem}

\author{Olena Karlova}

\ead{maslenizza.ua@gmail.com}

\address{Department of Mathematical Analysis, Faculty of Mathematics and Informatics, Chernivtsi National University, Kotsyubyns'koho str., 2, Chernivtsi, Ukraine}

\begin{abstract}
  We present $\sigma$-strongly functionally discrete mappings which expand the class of $\sigma$-discrete mappings and generalize Banach's theorem on analytically representable functions
\end{abstract}


\maketitle

\section{Introduction}

Let  $X$ be a topological space and $\mathcal A$ be a family of subsets of $X$. We define classes $\mathcal F_{\alpha}(\mathcal A)$ and $\mathcal G_{\alpha}(\mathcal A)$ in the following way: $\mathcal F_0(\mathcal A)=\mathcal A$, $\mathcal G_0(\mathcal A)=\{X\setminus A:A\in\mathcal A\}$ and for all $1\le\alpha<\omega_1$ let
 \begin{gather*}
   \mathcal F_{\alpha}(\mathcal A)=\Bigl\{\bigcap\limits_{n=1}^\infty A_n: A_n\in\bigcup\limits_{\beta<\alpha}\mathcal G_{\beta}(\mathcal A),\,\, n=1,2,\dots\Bigr\},\\
   \mathcal G_{\alpha}(\mathcal A)=\Bigl\{\bigcup\limits_{n=1}^\infty A_n: A_n\in\bigcup\limits_{\beta<\alpha}\mathcal F_{\beta}(\mathcal A),\,\, n=1,2,\dots\Bigr\}.
 \end{gather*}
If $\mathcal A$ is the collection of all /functionally/ closed subsets of $X$, then elements of $\mathcal F_{\alpha}(\mathcal A)$ or $\mathcal G_{\alpha}(\mathcal A)$ are called {\it sets of the $\alpha$'th /functionally/ multiplicative class} or  {\it sets of the $\alpha$'th /functionally/ additive class}, respectively; elements of the family $\mathcal F_{\alpha}(\mathcal A)\cap\mathcal G_{\alpha}(\mathcal A)$ are called {\it  /functionally/ ambiguous sets of the class $\alpha$}.

A mapping $f:X\to Y$ between topological spaces belongs to the {\it $\alpha$'th /functionally/ Lebesgue class}, if the preimage $f^{-1}(V)$ of any open set $V\subseteq Y$ is of the $\alpha$'th /functionally/ additive class $\alpha$ in $X$. The collection of all mappings of the $\alpha$'th /functionally/ Lebesgue class we denote by ${\rm H}_{\alpha}(X,Y)$ /${\rm K}_{\alpha}(X,Y)$/. Notice that ${\rm H}_\alpha(X,Y)={\rm K}_\alpha(X,Y)$ for any perfectly normal space  $X$ and any topological space $Y$.

By ${\rm C}(X,Y)$ we denote the class of all continuous mappings between $X$ and $Y$.

Let $\Phi_1(X,Y)={\rm H}_1(X,Y)$ and for all $1<\alpha<\omega_1$ the symbol $\Phi_{\alpha}(X,Y)$ stands for the collection of all mappings between $X$ and $Y$ which are pointwise limits of sequences of mappings from $\bigcup\limits_{\beta<\alpha}\Phi_{\beta}(X,Y)$. The next result is the classical Banach's theorem~\cite{BanachS:1931}.

\begin{thmm} Let $X$ be a metric space, $Y$ be a metric separable space and \mbox{$0<\alpha<\omega_1$}. Then
\begin{enumerate}
  \item[(i)] $\Phi_{\alpha}(X,Y)={\rm H}_{\alpha}(X,Y)$, if $\alpha<\omega_0$,

  \item[(ii)] $\Phi_{\alpha}(X,Y)={\rm H}_{\alpha+1}(X,Y)$, if $\alpha\ge \omega_0$.
\end{enumerate}
\end{thmm}

R.~Hansell in~\cite{Hansell:1971} introduced the concept of $\sigma$-discrete mapping as a convenient tool for the investigation of mappings with values in non-separable spaces. A mapping $f:X\to Y$  is $\sigma$-discrete if there exists a family  $\mathcal B=\bigcup\limits_{n=1}^\infty \mathcal B_n$ of subsets of $X$ such that every family $\mathcal B_n$ is discrete in $X$ and the preimage $f^{-1}(V)$ of any open set $V$ in $Y$ is a union of sets from $\mathcal B$. Class of all $\sigma$-discrete mappings between $X$ and $Y$ is denoted by $\Sigma(X,Y)$. In \cite[Theorem 7]{Hansell:1974} Hansell proved the following generalization of Banach's theorem.

\begin{thmm}\label{thm:Hansell} Let $X$ be a perfect space, $Y$ be a metric space and \mbox{$0<\alpha<\omega_1$}. Then
\begin{enumerate}
  \item[(i)] $\Phi_{\alpha}(X,Y)\cap\Sigma(X,Y)={\rm H}_{\alpha}(X,Y)\cap\Sigma(X,Y)$, if $\alpha<\omega_0$,

  \item[(ii)] $\Phi_{\alpha}(X,Y)\cap\Sigma(X,Y)={\rm H}_{\alpha+1}(X,Y)\cap\Sigma(X,Y)$, if $\alpha\ge \omega_0$.
\end{enumerate}
\end{thmm}

In this paper we develop technique from~\cite{Hansell:1974} and~\cite{Fos} and prove an analogue of Theorem~\ref{thm:Hansell} for an arbitrary topological space $X$. To do this we introduce classes of mappings $\Sigma_\alpha^f(X,Y)$. Namely, a  mapping $f:X\to Y$ belongs to $\Sigma_\alpha^f(X,Y)$, where $0\le\alpha<\omega_1$, if there exist a family $\mathcal B=\bigcup\limits_{n=1}^\infty \mathcal B_n$ of functionally  ambiguous sets of the class $\alpha$ in $X$ and a family $\mathcal U=\bigcup\limits_{n=1}^\infty \mathcal U_n$ of functionally open subsets of $X$, where $\mathcal U_n=(U_B:B\in\mathcal B_n)$, such that every family $\mathcal U_n$ is discrete in $X$, $\overline{B}\subseteq U_B$ for every $B\in\mathcal B_n$ and the preimage $f^{-1}(V)$ of any open set $V$ in $Y$ is a union of sets from $\mathcal B$. Properties of this class are studied in Section~\ref{sect:sigma}. Let us observe that the class $\Sigma_\alpha^f(X,Y)$ coincides with the class of all $\sigma$-discrete mappings of the  $\alpha$'th Lebesgue class in case $X$ is a perfectly normal space and $Y$ is a metric space. Auxiliary technical propositions are gathered in Section~\ref{sect:lemmas_sets}. The forth section contains three approximation lemmas which are crucial in the proof of the main theorem. In Section~\ref{sect:main_results} we present classes $\Lambda_\alpha$ which are close to classes $\Phi_\alpha$: let  $\Lambda_1(X,Y)=\Sigma_1^f(X,Y)$, and for all $1<\alpha<\omega_1$ let $\Lambda_{\alpha}(X,Y)$ be the collection of all mappings between $X$ and $Y$ which are pointwise limits of sequences of mappings from $\bigcup\limits_{\beta<\alpha}\Lambda_{\beta}(X,Y)$. The theorem below is the main result of the paper.
\begin{thmm} Let $X$ be a topological space, $Y$ be a metric space and \mbox{$0<\alpha<\omega_1$}. Then
\begin{enumerate}
  \item[(i)] $\Lambda_{\alpha}(X,Y)={\Sigma}_{\alpha}^f(X,Y)$, if  $\alpha<\omega_0$,

  \item[(ii)] $\Lambda_{\alpha}(X,Y)={\Sigma}_{\alpha+1}^f(X,Y)$, if $\alpha\ge \omega_0$.
\end{enumerate}
\end{thmm}
An example at the end of the fifth section shows that the assertion on  $X$ in Theorem~\ref{thm:Hansell} is essential.

\section{Properties of  $\sigma$-strongly functionally discrete mappings}\label{sect:sigma}

\begin{dfn}
 {\rm A family $\mathcal A=(A_i:i\in I)$ of subsets of a topological space $X$ is called
  \begin{enumerate}
    \item {\it discrete}, if every point of $X$ has an open neighborhood which intersects with at most one set from $\mathcal A$;

    \item {\it strongly discrete}, if there exists a discrete family $(U_i:i\in I)$ of open subsets of $X$ such that $\overline{A_i}\subseteq U_i$ for every $i\in I$;

    \item {\it strongly functionally discrete} or, briefly, {\it sfd-family}, if there exists a discrete family $(U_i:i\in I)$ of functionally open subsets of $X$ such that $\overline{A_i}\subseteq U_i$ for every $i\in I$;

    \item {\it well strongly functionally discrete} of {\it well sfd-family}, if there exist discrete families $(U_i:i\in I)$  of functionally open sets and $(F_i:i\in I)$ of functionally closed sets such that  $A_i\subseteq F_i\subseteq U_i$ for every $i\in I$.
  \end{enumerate}}
\end{dfn}

Notice that  (4) $\Rightarrow$ (3) $\Rightarrow$ (2) $\Rightarrow$ (1) for any $X$; a topological space $X$ is collectionwise normal if and only if every discrete family in $X$ is strongly discrete; if $X$ is normal then (2) $\Leftrightarrow$ (3).

\begin{dfn}
{\rm Let $\mathcal P$ be a property of a family of sets. A family $\mathcal A$ is called {\it a $\sigma$-$\mathcal P$-family} if  $\mathcal A=\bigcup\limits_{n=1}^\infty {\mathcal A}_n$, where every family $\mathcal A_n$ has the property $\mathcal P$.}
\end{dfn}

\begin{dfn}\label{def:base-for-a-function}
{\rm A family $\mathcal B$ of sets of a topological space $X$ is called {\it a base}  for a mapping $f:X\to Y$ if the preimage $f^{-1}(V)$ of an arbitrary open set  $V$ in $Y$ is a union of sets from $\mathcal B$.}
\end{dfn}

\begin{dfn}
{\rm  If a mapping $f:X\to Y$ has a base which is a $\sigma$-$\mathcal P$-family, then we say that $f$ is {\it a $\sigma$-$\mathcal P$ mapping}.}
\end{dfn}

The collection of all $\sigma$-$\mathcal P$ mappings between $X$ and $Y$ we will denote by
\begin{itemize}
  \item $\Sigma(X,Y)$, if $\mathcal P$ is a property of discreteness;

  \item $\Sigma^*(X,Y)$, if $\mathcal P$ is a property of a strong discreteness;

  \item $\Sigma^f(X,Y)$, if $\mathcal P$ is a property of a strong functional discreteness.
\end{itemize}

By $\Sigma^f_{\alpha}(X,Y)$ /$\Sigma^{f*}_{\alpha}(X,Y)$/ we denote the collection of all mappings between $X$ and $Y$ which has a $\sigma$-sfd base of functionally ambiguous /multiplicative/ sets of the class $\alpha$ in $X$.

\begin{rmk} For any spaces $X$ and $Y$ the following relations holds:
\begin{enumerate}
  \item   $\Sigma^{f*}_{\beta}(X,Y)\subseteq\Sigma^f_{\alpha}(X,Y)$m if $0\le\beta<\alpha<\omega_1$;

  \item $\Sigma^{f*}_{\beta}(X,Y)=\Sigma^f_{\beta+1}(X,Y)$m if $0\le \beta<\omega_1$;

  \item $\Sigma^{f}_0(X,Y)\subseteq C(X,Y)$.
\end{enumerate}

\end{rmk}

Let us observe that every continuous mapping $f:X\to Y$ is $\sigma$-strongly functionally discrete if either $X$, or $Y$ is a metric space, since every metric space has $\sigma$-sfd base of open sets. Clearly, every mapping with values in a second countable space is  $\sigma$-sfd. In~\cite{Hansell:1971} Hansell proved that any Borel measurable mapping $f:X\to Y$ is $\sigma$-discrete if $X$ is a complete metric space  and  $Y$ is a metric space.

\begin{lem} \label{union_SFD}
Let $0\le\alpha<\omega_1$, $X$ be a topological space, \mbox{$(U_i:i\in I)$} be a locally finite family of functionally open sets in $X$, $(B_i:i\in I)$ be a family of sets of the $\alpha$'th functionally additive /multiplicative/ class in  $X$ such that $B_i\subseteq U_i$ for every $i\in I$. Then the set $B=\bigcup\limits_{i\in I} B_i$ is of the $\alpha$'th functionally additive /multiplicative/ class~$\alpha$ in $X$.
\end{lem}

\begin{pf} For $\alpha=0$ we consider the case each $B_i$ is functionally closed and take a continuous function $f_i:X\to [0,1]$ such that  $B_i=f_i^{-1}(0)$ and $X\setminus U_i=f_i^{-1}(1)$, $i\in I$. Then the function $f(x)=\min\limits_{i\in I} f_i(x)$ is continuous and $B=f^{-1}(0)$.

Assume that our proposition is true for all $0\le\xi<\alpha$ and prove it for $\xi=\alpha$. If $\alpha$ is a limit ordinal then we take an increasing sequence of ordinals $(\alpha_n)_{n=1}^\infty$ which converges to $\alpha$. If $\alpha=\beta+1$ then we put $\alpha_n=\beta$ for every $n\in\mathbb N$.

Let $B_i$ be a set of the $\alpha$'th functionally additive class $\alpha$ and $(B_{i,n})_{n=1}^\infty$ be a sequence of sets of the $\alpha_n$'th functionally multiplicative classes such that $B_i=\bigcup\limits_{n=1}^\infty
B_{i,n}$.  By the inductive assumption the set $F_n=\bigcup\limits_{i\in I}
B_{i,n}$ belongs to  the $\alpha_n$'th functionally multiplicative class for every  $n$. Hence, the set $B=\bigcup\limits_{i\in I} B_i=\bigcup\limits_{n=1}^\infty F_n$ is of the $\alpha$'th functionally additive class.

Now assume that $B_i$ belongs to the $\alpha$'th functionally multiplicative class for every $i\in I$ and take a sequence $(B_{i,n})_{n=1}^\infty$ of sets of the $\alpha_n$'th functionally additive classes such that $B_i=\bigcap\limits_{n=1}^\infty
B_{i,n}$. Notice that each set $G_{i,n}=B_{i,n}\cap U_i$ is of the $\alpha_n$'th functionally additive class,
$G_{i,n}\subseteq U_i$ and $B_i=\bigcap\limits_{n=1}^\infty G_{i,n}$. Then $G_n=\bigcup\limits_{i\in I} G_{i,n}$ belongs to the $\alpha_n$'th functionally additive class for every  $n$. Hence, the set $B=\bigcap\limits_{n=1}^\infty G_n$ if of the $\alpha$'th functionally multiplicative class.
\end{pf}

\begin{cor}\label{cor:union_sfd}
For any  $0\le\alpha<\omega_1$ a union of an sfd-family of sets of the $\alpha$'th functionally additive /multiplicative/ class in a topological space is a set of the same class.
\end{cor}

\begin{lem}\label{lemma:special_base}
  Let $X$ be a topological space and $f\in\Sigma^f(X,Y)$. Then $f$ has a $\sigma$-sfd base $\mathcal B$ which is a union of a sequence of well  sfd-families.
\end{lem}

\begin{pf}
Let $\mathcal B'=\bigcup\limits_{n=1}^\infty \mathcal B_n'$ be a base for $f$, where $\mathcal B_n'$ is an sfd-family for every $n\in\mathbb N$. For all $n\in\mathbb N$ and $B\in\mathcal B_n'$ we take a functionally open set $U_{B,n}$ and a sequence of functionally closed sets  $(F_{B,m})_{m=1}^\infty$ such that the family $(U_{B,n}:B\in\mathcal B_n')$ is discrete, $B\subseteq U_{B,n}$ and $U_{B,n}=\bigcup\limits_{m=1}^\infty F_{B,m}$ for every $B\in\mathcal B'$. For all $n,m\in\mathbb N$ we put
  \begin{gather*}
  \mathcal B_{n,m}=(B\cap F_{B,m}:B\in\mathcal B_n')\quad\mbox{and}\quad \mathcal B=\bigcup\limits_{n,m=1}^\infty \mathcal B_{n,m}.
  \end{gather*}
It is easy to see that each $\mathcal B_{n,m}$ is well sfd-family and $\mathcal B$ is a base for $f$.
\end{pf}

\begin{thm}\label{LebesgueSFDisSFD0}
Let $0<\alpha<\omega_1$, $X$ be a topological space and $Y$  be a topological space with a $\sigma$-disjoint base. Then
  $$
  K_{\alpha}(X,Y)\cap\Sigma^f(X,Y)\subseteq\Sigma^f_{\alpha}(X,Y).
  $$
 \end{thm}

\begin{pf} Let $f\in K_{\alpha}(X,Y)\cap\Sigma^f(X,Y)$. According to Lemma~\ref{lemma:special_base} there exists a base $\mathcal B=\bigcup\limits_{m=1}^\infty\mathcal B_m$ for $f$ such that each $\mathcal B_m=(B_{i,m}:i\in I_m)$ is well sfd-family. For all $m$ and $i\in I_m$ we take a functionally open set $U_{i,m}$ and a functionally closed set $F_{i,m}$ in $X$ such that $B_{i,m}\subseteq F_{i,m}\subseteq U_{i,m}$ and the family  $(U_{i,m}:i\in I_m)$  is discrete.

Consider a $\sigma$-disjoint base $\mathcal V=\bigcup\limits_{n=1}^\infty\mathcal V_n$ of open sets in $Y$. Since $f\in K_{\alpha}(X,Y)$, for every  $V\in \mathcal V$ there exists a sequence $(A_{k,V})_{k=1}^\infty$ of sets of functionally multiplicative classes $<\alpha$ in $X$ such that
 $f^{-1}(V)=\bigcup\limits_{k=1}^\infty A_{k,V}$.
For  $m,n,k\in\mathbb N$ we put
$$
\mathcal B_{m,n,k}=(F_{i,m}\cap A_{k,V}:i\in I_m, V\in \mathcal V_n\,\,\,\mbox{and}\,\,\, B_{i,m}\subseteq f^{-1}(V)).
$$
Notice that each family $\mathcal B_{m,n,k}$ consists of functionally ambiguous sets of the class  $\alpha$  and is strongly functionally discrete in $X$, since the family $\mathcal B_m$ is strongly functionally discrete and for any nonempty set  $B_{i,m}\in \mathcal B_m$ there is at most one set $V\in\mathcal V_n$ such that $B_{i,m}\subseteq f^{-1}(V)$. Let
$$
\mathcal B_0=\bigcup\limits_{m,n,k=1}^\infty \mathcal B_{m,n,k}.
$$
We show that $\mathcal B_0$ is a base for $f$. Fix $V\in \mathcal V$ and verify that
$$
f^{-1}(V)=\bigcup\limits_{m,k=1}^\infty\bigcup\limits_{\mathop{i\in I_m}\limits_{B_{i,m}\subseteq f^{-1}(V)}}(F_{i,m}\cap A_{k,V}).
$$
Since $A_{k,V}\subseteq f^{-1}(V)$ for every $k$, the set in the right side of the equality is contained in $f^{-1}(V)$. On the other hand, if  $x\in f^{-1}(V)$ then $x\in A_{k,V}$ for some $k$. Moreover, $\mathcal B$ is a base for $f$, consequently, there are  $m$ and $i\in I_m$ such that  $x\in B_{i,m}\subseteq f^{-1}(V)$. Then $x\in  F_{i,m}$.
\end{pf}

\begin{prop}\label{Sigma0isK1}
Let $0<\alpha<\omega_1$, $X$ and $Y$ be topological spaces. Then
$$
\Sigma_{\alpha}^f(X,Y)\subseteq K_{\alpha}(X,Y)\cap\Sigma^f(X,Y).
$$
\end{prop}

\begin{pf}
Let  $f\in \Sigma_{\alpha}^f(X,Y)$. Clearly, $f\in \Sigma^f(X,Y)$. We show that $f\in K_{\alpha}(X,Y)$. Let $V$ be an open set in $Y$ and $\mathcal B=\bigcup\limits_{m=1}^\infty \mathcal B_m$ be a base for $f$ such that each family $\mathcal B_m$ is strongly functionally discrete in  $X$ and consists of functionally ambiguous sets of the class $\alpha$. Then there exists a subfamily $\mathcal B_V\subseteq\mathcal B$ such that $f^{-1}(V)=\bigcup \mathcal B_V$. For every $m\in\mathbb N$ we denote $\mathcal B_m'=(B\in\mathcal B_V:B\in\mathcal B_m)$. Corollary~\ref{cor:union_sfd} implies that every set $B_m=\bigcup\mathcal B_m'$ belongs to the $\alpha$'th functionally additive class in $X$. Moreover, $f^{-1}(V)=\bigcup\limits_{m=1}^\infty B_m$. Hence, $f\in K_{\alpha}(X,Y)$.
\end{pf}

Theorem~\ref{LebesgueSFDisSFD0} and Proposition~\ref{Sigma0isK1} imply
\begin{thm}\label{cor:sigma_F}
Let $0<\alpha<\omega_1$, $X$ be a topological space and $Y$  be a space with a $\sigma$-disjoint base.  Then
  $$
  K_{\alpha}(X,Y)\cap\Sigma^f(X,Y)=\Sigma^f_{\alpha}(X,Y).
  $$
\end{thm}

\begin{prop}\label{lemma:complex_map} Let $0\le\alpha<\omega_1$, $X$, $Y$ and $Z$ be topological spaces, $f\in\Sigma_{\alpha}^f(X,Y)$, $g\in\Sigma_{\alpha}^f(Y,Z)$ and let $h:X\to Y\times Z$ is defined by
$$
h(x)=(f(x),g(x))
$$
for every $x\in X$. Then $h\in\Sigma_{\alpha}^f(X,Y\times Z)$.
\end{prop}

\begin{pf}
  Let $\mathcal B_f=\bigcup\limits_{n=1}^\infty\mathcal B_{n,f}$ and $\mathcal B_g=\bigcup\limits_{m=1}^\infty\mathcal B_{m,g}$ be $\sigma$-sfd bases of functionally ambiguous sets of the class $\alpha$ for $f$ and $g$, respectively. For all $n,m\in\mathbb N$ we put
  $$
  \mathcal B_{n,m}=\{B_f\cap B_g: B_f\in\mathcal B_{n,f}, B_g\in\mathcal B_{m,g}\}.
 $$
It is easy to see that $\mathcal B=\bigcup\limits_{n,m=1}^\infty \mathcal B_{n,m}$ is a $\sigma$-sfd base for $h$ which consists of functionally ambiguous sets of the class  $\alpha$ in $X$.
\end{pf}

\begin{dfn}
{\rm We say that a family $(A_i:i\in I)$ is {\it a partition of a space $X$} if $X=\bigcup\limits_{i\in I}A_i$ and $A_i\cap A_j=\emptyset$ for $i\ne j$.}
\end{dfn}

\begin{prop}\label{lemma:restriction}
Let $0\le\alpha<\omega_1$, $(X_n:n\in\mathbb N)$ be a partition of a topological space $X$ by functionally ambiguous sets of the class $\alpha$, $(f_n)_{n=1}^\infty$ be a sequence of mappings from $\Sigma^f_\alpha(X,Y)$ and $f(x)=f_n(x)$ if $x\in X_n$ for some $n$. Then $f\in\Sigma_\alpha^f(X,Y)$.
\end{prop}

\begin{pf}
Let $\mathcal B_{n}$ be a $\sigma$-sfd base for a mapping $f_n$ which consists of functionally ambiguous sets of the class $\alpha$ in $X$. Let
  $$
  \mathcal B=\bigcup\limits_{n=1}^\infty(B\cap X_n: B\in\mathcal B_n).
  $$
It is easy to see that $\mathcal B$ is a $\sigma$-sfd base for $f$ which consists of functionally ambiguous sets of the $\alpha$'th class.
\end{pf}

\section{Auxiliary facts on functionally measurable sets}\label{sect:lemmas_sets}

The proofs of the next two lemmas are completely similar to the proofs of  Theorem 2 from \cite[p.~350]{Kuratowski:Top:1} and Theorem 2 from \cite[p.~357]{Kuratowski:Top:1}.

\begin{lem}\label{Lemma23} Let $0<\alpha<\omega_1$ and  $X$ be a topological space. Then for any disjoint sets $A,B\subseteq X$ of the $\alpha$'th functionally multiplicative class there exists a functionally ambiguous set $C$ of the class $\alpha$  such that $A\subseteq C\subseteq X\setminus B$.
\end{lem}

\begin{lem}\label{amb}
  If  $A$ is a functionally ambiguous set of the $(\alpha+1)$'th class in a topological space $X$, where $\alpha$ is a limit ordinal, then there exists a sequence $(A_n)_{n=1}^\infty$ of functionally ambiguous sets of classes $<\alpha$ such that
  \begin{gather}
    A=\bigcup\limits_{n=1}^\infty \bigcap\limits_{k=0}^\infty A_{n+k}=\bigcap\limits_{n=1}^\infty \bigcup\limits_{k=0}^\infty A_{n+k}.
  \end{gather}
\end{lem}

The definition of sfd-family easily implies the following fact.
\begin{lem}\label{lemma:technical-1}
Let $\mathcal A_1$,\dots, $\mathcal A_n$ be sfd-families of subsets of a topological space $X$,   $A_k=\bigcup\mathcal A_k$ for $k=1,\dots,n$ and the family $(A_k:k=1,\dots,n)$ is strongly functionally discrete. Then the family $\mathcal A=\bigcup\limits_{k=1}^n\mathcal A_k$ is strongly functionally discrete.
\end{lem}

\begin{lem}\label{lemma:technical-2}
Let $0<\alpha<\omega_1$,  $\mathcal A$ be a disjoint $\sigma$-sfd family of functionally additive sets of the $\alpha$'th class in a topological space $X$. Then for any $A\in\mathcal A$ there exists an increasing sequence $(D_{n}^A)_{n=1}^\infty$ of functionally ambiguous sets of the class $\alpha$ such that $A=\bigcup\limits_{n=1}^\infty D_{n}^A$ and the family $(D_{n}^A:A\in\mathcal A)$ is strongly functionally discrete for every  $n\in\mathbb N$.

If $\alpha=\beta+1$, then every set $D_n^A$ can be chosen from the $\beta$'th functionally multiplicative class.
\end{lem}

\begin{pf} Let $\alpha=1$ and $\mathcal A=\bigcup\limits_{k=1}^\infty\mathcal A_k$, where $\mathcal A_k$ is an sfd-family of sets of the first functionally additive class and $\bigl(\cup\mathcal A_k\bigr)\cap \bigl(\cup\mathcal A_j\bigr)=\emptyset$ for $k\ne j$. For every $A\in\mathcal A$ we take an increasing sequence $(B_{n}^A)_{n=1}^\infty$ of functionally closed sets such that $A=\bigcup\limits_{n=1}^\infty B_{n}^A$. Now for all $A\in\mathcal A$ and $n\in\mathbb N$ we put
  $$
  F_{n}^A=\left\{\begin{array}{ll}
                   B_{n}^A, &  \mbox{if\,\,} A\in\mathcal A_k\mbox{\,\,\,for\,\,\,}k\le n, \\
                   \emptyset, & \mbox{if\,\,} A\in\mathcal A_k\mbox{\,\,\,for\,\,\,}k> n.
                 \end{array}
  \right.
  $$
  Then $(F_{n}^A:A\in\mathcal A)=\bigcup\limits_{k\le n}(B_{n}^A:A\in\mathcal A_k)$ for every $n$. Since every family $\mathcal B_k=(B_{n}^A:A\in\mathcal A_k)$ is strongly functionally discrete, the set $B_k=\bigcup\limits\mathcal B_k$ is functionally closed. Moreover, $B_k\cap B_m=\emptyset$ for all $k\ne m$. Since $(B_k:k=1,\dots,n)$ is an sfd-family, Lemma~\ref{lemma:technical-1} implies that $(F_{n}^A:A\in\mathcal A)$ is also sfd-family. Moreover, $A=\bigcup\limits_{n=1}^\infty F_{n}^A$ for every $A\in\mathcal A$.

 Assume that the assertion of lemma is true for all $1\le\xi<\alpha$ and verify it for $\xi=\alpha$. Consider a disjoint sequence of sfd-families $\mathcal A_k$ which consist of sets of the $\alpha$'th functionally additive class. For every $A\in\mathcal A$ we take an increasing sequence $(B_n^A)_{n=1}^\infty$ such that $A=\bigcup\limits_{n=1}^\infty B_{n}^A$. We may assume that every set $B_n^A$ belongs to the $\alpha_n$'th functionally multiplicative class, where  $\alpha_n=\beta$ for every  $n\in\mathbb N$ if $\alpha=\beta+1$,  and $(\alpha_n)_{n=1}^\infty$ is an increasing sequence of ordinals such that $\alpha={\rm sup}\,\alpha_n$ if $\alpha$ is a limit ordinal.

  Fix $n\in\mathbb N$ and for every $k=1,\dots,n$ we denote $\mathcal B_{k,n}=(B_n^A:A\in\mathcal A_k)$ and $B_{k,n}=\cup\mathcal B_{k,n}$. Since $B_{1,n},\dots,B_{n,n}$ are mutually disjoint sets of the $\alpha_n$'th functionally multiplicative class, Lemma~\ref{Lemma23} implies that there exist mutually disjoint functionally ambiguous sets $C_{1,n},\dots,C_{n,n}$ of the class $\alpha_n$ such that $B_{k,n}\subseteq C_{k,n}$ for every $k=1,\dots,n$. By the inductive assumption $C_{k,n}=\bigcup\limits_{m=1}^\infty C_{k,n,m}$ and for every $m\in\mathbb N$ the family $(C_{k,n,m}:k=1,\dots,n)$ is strongly functionally discrete and consists of functionally ambiguous sets of the class $\alpha_n$.

Now for all $n,m\in\mathbb N$ and $A\in\mathcal A$ we put
  $$
  D_{n,m}^A=\left\{\begin{array}{ll}
                     B_n^A\cap C_{k,m,n}, & \mbox{if\,\,} A\in\mathcal A_k\mbox{\,\,\,for\,\,\,}k\le n, \\
                     \emptyset, & \mbox{if\,\,} A\in\mathcal A_k\mbox{\,\,\,for\,\,\,}k> n.
                   \end{array}
  \right.
  $$
Then $(D_{n,m}^A:A\in\mathcal A)=\bigcup\limits_{k\le n}(B_n^A\cap C_{k,m,n}:A\in\mathcal A_k)$.

Fix $n,m\in\mathbb N$ and for every  $k=1,\dots,n$ we put
  $$
  \mathcal D_k=(B_n^A\cap C_{k,m,n}:A\in\mathcal A_k).
  $$
Notice that every family $\mathcal D_k$ is strongly functionally discrete and consists of functionally ambiguous sets of the class  $\alpha$. Then the set $D_k=\bigcup\limits\mathcal D_k$ is functionally ambiguous of the class $\alpha$ for $k=1,\dots,n$. Moreover, $D_k\cap D_m=\emptyset$ for all $k\ne m$. Since the family $(D_k:k=1,\dots,n)$ is strongly functionally discrete, Lemma~\ref{lemma:technical-1} implies that $(D_{n,m}^A:A\in\mathcal A)$ is an sfd-family. Moreover, $A=\bigcup\limits_{n,m=1}^\infty D_{n,m}^A$ for every $A\in\mathcal A$.

In case $\alpha=\beta+1$ for all $n,m\in\mathbb N$ and $A\in\mathcal A$ we choose an increasing sequence  $(D_{n,m,k}^A)_{k=1}^\infty$ of functionally multiplicative class $\beta$ such that $D_{n,m}^A=\bigcup\limits_{k=1}^\infty D_{n,m,k}^A$. Clearly, $(D_{n,m,k}^A:A\in\mathcal A)$ is an sfd-family for all $n,m,k$ and $A=\bigcup\limits_{n,m,k=1}^\infty D_{n,m,k}^A$.
 \end{pf}

\begin{lem}\label{lemma:technical-3}
  Let $0\le \alpha<\omega_1$, $X$ be a topological space, $\mathcal A$ be an $\sigma$-sfd family of sets of the $\alpha$'th functionally multiplicative class such that $\bigcup\mathcal A=X$. Then there exists a sequence $(\mathcal A_n)_{n=1}^\infty$ of families of sets of the $\alpha$'th functionally multiplicative class such that
  \begin{enumerate}
    \item $\bigcup\limits_{n=1}^\infty \mathcal A_n\prec\mathcal A$,

    \item $\mathcal A_{n}\prec\mathcal A_{n+1}$,

    \item $\mathcal A_n$ is an sfd-family,

    \item $\bigcup\bigcup\limits_{n=1}^\infty \mathcal A_n=X$
  \end{enumerate}
  for every $n\in\mathbb N$.
\end{lem}

\begin{pf}
  Let $\mathcal A=\bigcup\limits_{k=1}^\infty \mathcal B_k$ and let $\mathcal B_k$ be an sfd-family of sets of the $\alpha$'th functionally multiplicative class. For every $k\in\mathbb N$ we put
  $$
  \mathcal C_k=(B\setminus \bigcup\limits_{j<k}\bigcup \mathcal B_j: B\in\mathcal B_k)\quad\mbox{and}\quad \mathcal C=\bigcup\limits_{k=1}^\infty \mathcal C_k.
  $$
  Then $\mathcal C$ is a disjoint $\sigma$-sfd family of sets of the $(\alpha+1)$'th functionally additive class, $\mathcal C\prec \mathcal A$ and $\bigcup\mathcal C=X$. By Lemma~\ref{lemma:technical-2} for every $C\in\mathcal C$ there exists an increasing sequence $(D_{n}^C)_{n=1}^\infty$ of sets of the $\alpha$'th functionally multiplicative class such that  $C=\bigcup\limits_{n=1}^\infty D_{n}^C$ and the family $(D_{n}^C:C\in\mathcal C)$ is strongly functionally discrete for every $n\in\mathbb N$.

It remains to put $\mathcal A_n=(\bigcup\limits_{k\le n}D_{k}^C:C\in\mathcal C)$ for $n\in\mathbb N$.
\end{pf}

\section{Approximation lemmas}

\begin{dfn}
We say that a sequence $(f_n)_{n=1}^\infty$ of mappings $f_n:X\to Y$ {\it is stably convergent to $f:X\to Y$} and denote this fact by $f_n\stackrel{{\rm st}}{\mathop{\longrightarrow}} f$, if for every $x\in X$ there exists $n_0\in\mathbb N$ such that $f_n(x)=f(x)$ for all $n\ge n_0$.
\end{dfn}

If $A\subseteq Y^X$, then the symbol $\overline{A}^{\,\,{\rm st}}$ stands for the set of all stable limits of sequences from $A$.

\begin{lem}\label{lemma:approx:limit}
Let $X$, $Y$ be topological spaces, $((B_{i,m}:i\in I_m))_{m=1}^\infty$ be a sequence of sfd-families of sets of the $(\alpha+1)$'th functionally ambiguous sets in $X$, $\alpha$ be a limit ordinal, let the family $(B_m=\bigcup\limits_{i\in I_m}B_{i,m}:m\in\mathbb N)$ be a partition of $X$, $((y_{i,m}:i\in I_m))_{m=1}^\infty$ be a sequence of points from $Y$ and let $f:X\to Y$ is defined by
$$
f(x)=y_{i,m},
$$
if $x\in B_{i,m}$ for some $m\in\mathbb N$ and $i\in I_m$. Then $f\in \overline{\Sigma_{<\alpha}^f(X,Y)}^{\,\,{\rm st}}$.
\end{lem}

\begin{pf}
Fix $m\in\mathbb N$. Since $B_m$ is functionally ambiguous set of the class $(\alpha+1)$ in  $X$ by Corollary~\ref{union_SFD}, Lemma~\ref{amb} implies that there exists a sequence $(C_{m,n})_{n=1}^\infty$ of functionally ambiguous sets of classes $<\alpha$ such that
 \begin{gather}\label{eq:amb}
  B_m=\bigcup\limits_{n=1}^\infty \bigcap\limits_{k=0}^\infty C_{m,n+k}=\bigcap\limits_{n=1}^\infty \bigcup\limits_{k=0}^\infty C_{m,n+k}.
 \end{gather}
Moreover, there exists a discrete family $(U_{i,m}:i\in I_m)$ of functionally open sets in $X$ such that $B_{i,m}\subseteq U_{i,m}$ for every  $i\in I_m$.

For all $m,n\in\mathbb N$ we put
 $$
 D_{m,n}=C_{m,n}\setminus\bigcup\limits_{k<m} C_{k,n}.
 $$
Notice that every $D_{m,n}$ is a functionally ambiguous set of a class $<\alpha$. Moreover,
 \begin{gather}\label{eq:cond_on_Dmn}
   \bigl(\forall x\in X\bigr)\, \bigl(\exists m_x\in\mathbb N\bigr)\, \bigl(\exists n_{x}\in\mathbb N\bigr)\, \bigl(\forall n\ge n_{x}) \, \bigl(x\in D_{m_x,n}).
 \end{gather}
Indeed, if $x\in X$, then there exists a unique number $m_x$ such that $x\in B_{m_x}$ and $x\not\in B_k$ for all $k\ne {m_x}$. Then the equality~\ref{eq:cond_on_Dmn}) implies that there are numbers $N_1,\dots,N_{m_x}$ such that
 $$
 x\not\in \bigcup\limits_{n\ge N_k} C_{k,n}\,\,\mbox{if}\,\, k<m_x\,\,\,\,\mbox{and}\,\,\,\,  x\in \bigcap\limits_{n\ge N_{m_x}} C_{m_x,n}.
 $$
Hence, for all $n\ge n_{x}=\max\{N_1,\dots,N_{m_x}\}$ we have $x\in D_{m_x,n}$.

Let $y_0$ be any point from $\{y_{i,1}:i\in I_1\}$. Fix $n\in\mathbb N$ and for all $x\in X$ let
 $$
 f_n(x)=\left\{\begin{array}{ll}
                 y_{i,m}, & \mbox{if}\,\,\, x\in D_{m,n}\cap U_{i,m}\,\,\,\mbox{for some}\,\,\, m<n \,\,\,\mbox{and}\,\,\, i\in I_m,\\
                 y_0, & \mbox{otherwise}.
               \end{array}
 \right.
 $$
Observe that $f_n:X\to Y$ is defined correctly, since the family $(U_{i,m}:i\in I_m)$ is discrete for every $m$ and the family $(D_{m,n}:m<n)$ if disjoint.

We show that $f_n\in \Sigma_{<\alpha}^f(X,Y)$. For $m=1,\dots,n-1$ let $\mathcal B_m=(D_{m,n}\cap U_{i,m}: i\in I_m)$ and let $\mathcal B_n$ be a family which consists of the set  $\bigl(X\setminus \bigcup\bigcup\limits_{m<n}\mathcal B_m\bigr)$. Clearly, $\mathcal B=\bigcup\limits_{m=1}^n\mathcal B_m$ is $\sigma$-sfd family of functionally ambiguous sets of classes $<\alpha$. It follows from the definition of $f_n$ that $\mathcal B$ is a base for $f_n$.

It remains to prove that $f_n\stackrel{{\rm st}}{\mathop{\longrightarrow}} f$ on $X$. Indeed, if $x\in X$, then there exist $m\in\mathbb N$ and $i\in I_m$ such that $x\in B_{i,m}\subseteq U_{i,m}$. Then $f(x)=y_{i,m}$. It follows from~(\ref{eq:cond_on_Dmn}) that there exists a number $n_0> m$ with $x\in D_{m,n}$ for all $n\ge n_{0}$. Then $f_n(x)=y_{i,m}$ for all $n\ge n_0$. Hence, $f_n(x)=f(x)$ for all $n\ge n_0$.
\end{pf}

\begin{lem}\label{lemma:uniform-high}
Let $\alpha<\omega_1$ be a limit ordinal, $X$ be a topological space, $Y$ be a metric space and $f\in\Sigma_{\alpha+1}^{f}(X,Y)$. Then there exists a sequence of mappings  $f_n\in\overline{\Sigma_{<\alpha}^f(X,Y)}^{\,\,\rm st}$ which is uniformly convergent to $f$ on $X$.
\end{lem}

\begin{pf}
   Let $\mathcal B$ be a $\sigma$-sfd base for $f$ which consists of functionally ambiguous sets of the class $(\alpha+1)$ in $X$. For every  $n\in\mathbb N$ we consider a covering $\mathcal U_n$ of $Y$ by open balls of diameters $<\frac{1}{n}$ and put
\begin{gather*}
\mathcal B_n=(B\in\mathcal B: \exists U\in \mathcal U_n \,\,\,|\,\,\, B\subseteq f^{-1}(U)).
\end{gather*}
Then $\mathcal B_n$ is a $\sigma$-sfd family of functionally ambiguous sets of the class $(\alpha+1)$, ${\rm diam\,\,}f(B)<\frac{1}{n}$ for every  $B\in\mathcal B_n$ and $X=\cup \mathcal B_n$ for every $n\in\mathbb N$.

Fix $n\in\mathbb N$ and let $\mathcal B_n=\bigcup\limits_{m=1}^\infty \mathcal B_{n,m}$, where $\mathcal B_{n,m}$ is an sfd-family of functionally ambiguous sets of the class  $(\alpha+1)$. We put $\mathcal A_{n,1}=\mathcal B_{n,1}$ and $\mathcal A_{n,m}=\bigl(B\setminus\bigl(\bigcup\bigcup\limits_{k<m}\mathcal B_{n,k}\bigr):B\in\mathcal B_{n,m}\bigr)$ for $m>1$. Notice that for every $m\in\mathbb N$ the set $A_{n,m}=\cup\mathcal A_{n,m}$ is functionally ambiguous of the class  $(\alpha+1)$ and the family $(A_m:m\in\mathbb N)$ is a partition of $X$. For every $A\in\mathcal A_{n,m}$ we choose an arbitrary point $y_{n,m}^A\in f(A)$. We define a mapping $f_n:X\to Y$  by
$$
f_n(x)=y_{n,m}^A,
$$
if $x\in A$ for some $m\in\mathbb N$ and $A\in\mathcal A_{n,m}$. Then $f_n\in\overline{\Sigma_{<\alpha}^f(X,Y)}^{\,\,\rm st}$ by Lemma~\ref{lemma:approx:limit}.

It remains to verify that $(f_n)_{n=1}^\infty$ converges uniformly to $f$. Indeed, if $x\in X$ and $n\in\mathbb N$, then $f_n(x)=y_{n,m}^A\in f(A)$ for some
 $m\in\mathbb N$ and $A\in\mathcal A_{n,m}$. Since $A\subseteq B$ for some $B\in\mathcal B_{n,m}$,
$$
d(f(x),f_n(x))\le{\rm diam} f(B)<\frac{1}{n},
$$
which completes the proof.
\end{pf}

\begin{lem}\label{prop:main_technical:finite}
Let $0<\alpha<\omega_1$, $X$ be a topological space, $(Y,d)$ be a metric space, \mbox{$f:X\to Y$} be a mapping,
${\mathcal A}_1,\dots,{\mathcal A}_n$ be families of subsets of $X$ such that 
\begin{enumerate}
  \item[(i)] $\mathcal A_k$ is an sfd-family of sets of the $\alpha$'th functionally multiplicative class;

  \item[(ii)] ${\mathcal A}_{k+1}\prec {\mathcal A}_k$ for $k<n$;

  \item[(iii)] ${\rm diam}(f(A))<\frac{1}{2^{k+2}}$ for all $A\in {\mathcal A}_k$
\end{enumerate}
for every $k=1,\dots,n$.  Then there exists a mapping $g\in \Sigma^f_{\alpha}(X,Y)$ such that the inclusion $x\in \cup {\mathcal A}_k$ for   $k=1,\dots,n$ implies the inequality 
\begin{gather}\label{prop:main_technical_ineq}
d(f(x),g(x))<\frac{1}{2^k}.
\end{gather}
\end{lem}

\begin{pf} Let $\mathcal A_k=(A_{i,k}:i\in I_k)$ and $(U_{i,k}:i\in I_k)$ be discrete families of functionally open sets in  $X$ such that  $A_{i,k}\subseteq U_{i,k}$ for every $i\in I_k$ and $k=1,\dots,n$.

By Lemma~\ref{Lemma23} there exists a family $(B_{i,1}:i\in I_1)$ of functionally ambiguous sets of the class  $\alpha$ such that $A_{i,1}\subseteq B_{i,1}\subseteq U_{i,1}$. Since $\mathcal A_2\prec \mathcal A_1$, for every $i\in I_2$ there exists a unique $j\in I_1$ such that $A_{i,2}\subseteq A_{j,1}$. Notice that $U_{i,2}\cap B_{j,1}$ is a functionally ambiguous set of the class $\alpha$ and applying Lemma~\ref{Lemma23} we obtain a functionally ambiguous set  $B_{i,2}$ of the class  $\alpha$ such that $A_{i,2}\subseteq B_{i,2}\subseteq U_{i,2}\cap B_{j,1}$. Proceeding in this way we obtain a sequence $((B_{i,k}:  i\in I_k))_{k=1}^n$ of families of subsets of  $X$ such that 
\begin{itemize}
\item $B_{i,k}\subseteq U_{i,k}$ for every $i\in I_k$;

  \item $B_{i,k}$ is a functionally ambiguous set of the class $\alpha$ for every $i\in I_k$;

  \item for every $k<n$ and for every $i\in I_{k+1}$  there exists a unique $j\in I_k$ such that 
  \begin{equation}\label{prop:main_technical_cond1}
  A_{i,k+1}\subseteq A_{j,k},
  \end{equation}
  \begin{equation}\label{prop:main_technical_cond2}
   A_{i,k+1}\subseteq B_{i,k+1}\subseteq B_{j,k}.
  \end{equation}
\end{itemize}
for all $k=1,\dots,n$. 
Observe that for every  $k$ the set 
$$
B_k=\bigcup\limits_{i\in I_k}B_{i,k}
$$
is functionally ambiguous of the class  $\alpha$ according to Corollary~\ref{union_SFD}.

We take any points $y_0\in f(X)$ and $y_{i,k}\in f(A_{i,k})$  for every $k$ and $i\in I_k$. For all $x\in X$ we put
$$
g_0(x)=y_0.
$$
Assume that for some $k<n$ we have already defined mappings  $g_1,\dots,g_k$ from  $\Sigma_\alpha^f(X,Y)$ such that 
\begin{gather}\label{c4}
g_k(x)=\left\{\begin{array}{ll}
                g_{k-1}(x), & \mbox{if\,\,\,}x\in X\setminus B_k, \\
                y_{i,k}, & \mbox{if\,\,\,}x\in B_{i,k}\mbox{\,\,for some\,\,} i\in I_k.
              \end{array}
\right.
\end{gather}
We put
$$
g_{k+1}(x)=\left\{\begin{array}{ll}
                g_{k}(x), & \mbox{if\,\,\,}x\in X\setminus B_{k+1}, \\
                y_{i,k+1}, & \mbox{if\,\,\,}x\in B_{i,k+1}\mbox{\,\,for some\,\,} i\in I_{k+1}.
              \end{array}
\right.
$$
Then $g_{k+1}\in\Sigma_\alpha^f(X,Y)$ by Lemma~\ref{lemma:restriction}. Repeating inductively this process we obtain mappings  $g_1,\dots,g_n$ from $\Sigma_\alpha^f(X,Y)$ each of which satisfies~(\ref{c4}).

Now we prove that 
\begin{gather}\label{c3}
  d(g_{k+1}(x),g_k(x))<\frac{1}{2^{k+2}}
\end{gather}
for all $0\le k<n$ and $x\in X$. Indeed, if  $x\in X\setminus B_{k+1}$, then  $g_{k+1}(x)=g_k(x)$ and \mbox{$d(g_{k+1}(x),g_k(x))=0$}. Assume that  $x\in B_{i,k+1}$ for some $i\in I_{k+1}$ and choose  $j\in I_k$ such that (\ref{prop:main_technical_cond1}) and~(\ref{prop:main_technical_cond2}) holds. Then $g_{k+1}(x)=y_{i,k+1}$ and $g_k(x)=y_{j,k}$. Since $f(A_{i,k+1})\subseteq f(A_{j,k})$,  $y_{i,k+1}\in f(A_{j,k})$. Hence, $d(g_{k+1}(x), g_k(x))\le {\rm diam}(f(A_{j,k}))<\frac{1}{2^{k+2}}$.

We put $g=g_n$ and show that~(\ref{prop:main_technical_ineq}) holds. Let $1\le k\le n$ and $x\in \cup\mathcal A_k$. Then $x\in A_{i,k}\subseteq B_{i,k}$ for some $i\in I_k$. It follows that $g_k(x)=y_{i,k}$ and consequently 
$$
d(f(x),g_k(x))\le {\rm diam}(f(A_{i,k}))<\frac{1}{2^{k+2}}.
$$
Taking into account~(\ref{c3}) we obtain that 
\begin{gather*}
  d(f(x),g_n(x))\le d(f(x),g_k(x))+\sum\limits_{i=k}^{n-1} d(g_i(x),g_{i+1}(x))<\frac{1}{2^{k+2}}+\frac{1}{2^{k+1}}<\frac{1}{2^k}.
\end{gather*}
\end{pf}

\section{A generalization of Banach's theorem}\label{sect:main_results}

Let $X$, $Y$ be topological spaces and $A\subseteq Y^X$. By $\overline{A}^{\,\,{\rm p}}$ we denote the set of all pointwise limits of sequences of mappings from $A$.

We put
$$
\Lambda_1(X,Y)=\Sigma_1^f(X,Y)
$$
and for all $1<\alpha<\omega_1$ let
$$
\Lambda_{\alpha}(X,Y)=\overline{\bigcup\limits_{\beta<\alpha}\Lambda_{\beta}(X,Y)}^{\,\,{\rm p}}.
$$
Clearly, 
$$
\Lambda_{\beta}(X,Y)\subseteq \Lambda_{\alpha}(X,Y),
$$
if $\beta\le\alpha$. Moreover,
$$
\Lambda_{\alpha+1}(X,Y)=\overline{\Lambda_{\alpha}(X,Y)}^{\,\,{\rm p}},
$$
and if $\alpha={\rm sup\,}\alpha_n$ is a limit ordinal, then
$$
\Lambda_{\alpha}(X,Y)=\overline{\bigcup\limits_{n=1}^\infty\Lambda_{\alpha_n}(X,Y)}^{\,\,{\rm p}}.
$$

\begin{thm}\label{finite_classes}
Let $X$ be a topological space, $(Y,d)$ be a metric space. Then 
 \begin{enumerate}
   \item[(i)] $\Sigma_{\alpha}^f(X,Y)\subseteq \Lambda_{\alpha}(X,Y)$, if $1\le \alpha\le \omega_0$;

   \item[(ii)] $\Sigma_{\alpha+1}^f(X,Y)\subseteq \Lambda_{\alpha}(X,Y)$, if $\omega_0\le \alpha<\omega_1$.
 \end{enumerate}
\end{thm}

\begin{pf} The proposition is obvious for $\alpha=1$. Assume it is true for all $1\le \beta<\alpha$ and prove the proposition for $\beta=\alpha$.

Let $\alpha<\omega_0$ and let $f$ be a mapping from $\Sigma_{\alpha}^{f}(X,Y)=\Sigma_{\alpha-1}^{f*}(X,Y)$ with a $\sigma$-sfd base $\mathcal B$ which consists of sets of the $(\alpha-1)$'th functionally multiplicative class in  $X$. For every $k\in\mathbb N$ we consider a covering $\mathcal U_k$ of $Y$ by open balls of diameters $<\frac{1}{2^k}$ and put
$$
\mathcal B_k=(B\in\mathcal B: \exists U\in \mathcal U_k \,\,\,|\,\,\, B\subseteq f^{-1}(U)).
$$
Then $\mathcal B_k$ is a $\sigma$-sfd family and $X=\cup \mathcal B_k$ for every $k$. By Lemma~\ref{lemma:technical-3} for every $k\in\mathbb N$ there exists a sequence  $(\mathcal B_{k,n})_{n=1}^\infty$ of  sfd families of sets of the $(\alpha-1)$'th functionally multiplicative class in $X$ such that $\bigcup\limits_{n=1}^\infty \mathcal B_{k,n}\prec\mathcal \mathcal B_k$, $\mathcal B_{k,n}\prec\mathcal B_{k,n+1}$ and $\bigcup\bigcup\limits_{n=1}^\infty \mathcal B_{k,n}=X$ for every $n\in\mathbb N$. For all $k,n\in\mathbb N$ we put
$$
\mathcal F_{k,n}=(B_1\cap\dots\cap B_k: B_m\in \mathcal B_{m,n}, 1\le m\le k).
$$
Notice that every family $\mathcal F_{k,n}$ is strongly functionally discrete, consists of sets of the $(\alpha-1)$'th functionally multiplicative class and
\begin{enumerate}
  \item[(a)] $\mathcal F_{k+1,n}\prec \mathcal F_{k,n}$,

  \item[(b)] $\mathcal F_{k,n}\prec \mathcal F_{k,n+1}$,

  \item[(c)]  $\bigcup\limits_{n=1}^\infty \mathcal F_{k,n}=X$.
\end{enumerate}
For every $n\in\mathbb N$ we apply Lemma~\ref{prop:main_technical:finite} to the mapping $f$ and to the families $\mathcal F_{1,n}$, $\mathcal F_{2,n}$,\dots,$\mathcal F_{n,n}$ and obtain a sequence of mappings $g_n\in \Sigma^f_{\alpha-1}(X,Y)$ such that the inclusion $x\in \cup\mathcal F_{k,n}$ for  $k\le n$ implies the inequality
$$
d(f(x),g_n(x))<\frac{1}{2^k}.
$$

We prove that $g_n\to f$ pointwise on $X$. Fix  $x\in X$ and $\varepsilon>0$. Let $k\in\mathbb N$ be a number such that $\frac{1}{2^k}<\varepsilon$.  Conditions~(b) and (c) imply that there exists $n_0\ge k$ such that $x\in \mathcal F_{k,n}$ for all $n\ge n_0$. Then for all $n\ge n_0$ we have
$d(f(x),g_n(x))<\frac{1}{2^k}<\varepsilon$.

Since $g_n\in \Lambda_{\alpha-1}(X,Y)$ for every $n$ by the inductive assumption, $f\in \Lambda_{\alpha}(X,Y)$. Therefore, the proposition~(i) is proved for all  $\alpha<\omega_0$.

Now let  $\alpha\ge\omega_0$ be a limit ordinal and $f\in\Sigma_{\alpha+1}^f(X,Y)$. Lemma~\ref{lemma:uniform-high} implies that there exists a sequence of mappings  $f_n\in\overline{\Sigma_{<\alpha}^f(X,Y)}^{\,\,\rm st}$ which converges uniformly to $f$ on $X$. Without loss of generality we may assume that 
$$
d(f_{n+1}(x),f_n(x))<\frac{1}{2^n}
$$
for all $x\in X$ and $n\in\mathbb N$. For every $n\in\mathbb N$ there exists a sequence of mappings $f_{n,m}\in\Sigma_{<\alpha}^f(X,Y)$ such that 
$$
f_{n,m}\mathop{\stackrel{{\rm st}}{\mathop{\longrightarrow}}}\limits_{m\to\infty} f_n
$$
on $X$. For every $x\in X$ we put $h_{0,m}(x)=f_{0,m}(x)$ for every $m\in\mathbb N$. Suppose that we have already defined sequences of mappings $(h_{k,m})_{m=1}^\infty$ for some $n\in\mathbb N$ and for every $0\le k\le n$ such that
\begin{enumerate}
  \item[(a)] $h_{k,m}\stackrel{{\rm st}}{\mathop{\longrightarrow}} f_k$ for all  $0\le k\le n$;

  \item[(b)] $h_{k,m}\in\Sigma^f_{<\alpha}(X,Y)$ for all $0\le k\le n$ and $m\in\mathbb N$;

  \item[(c)] $d(h_{k+1,m}(x),h_{k,m}(x))<\frac{1}{2^{k}}$ for all $0\le k<n$, $m\in\mathbb N$ and $x\in X$.
\end{enumerate}
For every  $m\in\mathbb N$ we consider the set
$$
A_m=\left\{x\in X: d(f_{n+1,m}(x),h_{n,m}(x))<\frac{1}{2^{n}}\right\}
$$
and put
$$
h_{n+1,m}(x)=\left\{\begin{array}{ll}
                      f_{n+1,m}(x), & \mbox{if}\,\,\, x\in A_m, \\
                      h_{n,m}(x), & \mbox{if}\,\,\, x\not\in A_m.
                    \end{array}
\right.
$$

We check the condition (a) for the sequence $(h_{n+1,m})_{m=1}^\infty$. If $x\in X$, then there exists a number $m_0$ such that $f_{n+1,m}(x)=f_{n+1}(x)$ and $f_{n,m}(x)=f_n(x)$ for all $m\ge m_0$. Then
$$
d(f_{n+1,m}(x),h_{n,m}(x))=d(f_{n+1}(x),f_n(x))<\frac{1}{2^{n}}
$$
for all $m\ge m_0$. Hence, $x\in A_m$ for all $m\ge m_0$. Consequently, $h_{n+1,m}(x)=f_{n+1,m}(x)$, which implies that 
$$
h_{n+1,m}\mathop{\stackrel{{\rm st}}{\mathop{\longrightarrow}}}\limits_{m\to\infty} f_{n+1}
$$
on $X$.

Now we check the condition (b). For every $m\in\mathbb N$ the mapping $\varphi(x)=d(f_{n+1,m}(x),h_{n,m}(x))$ belongs to a class $\Sigma_{<\alpha}^f(X,\mathbb R)$ according to Proposition~\ref{lemma:complex_map}. Hence, every set $A_m=\varphi^{-1}((-\infty,\frac{1}{2^{n}}))$ is functionally ambiguous of a class $<\alpha$. Thus, $h_{n+1,m}\in\Sigma^f_{<\alpha}(X,Y)$ by Proposition~\ref{lemma:restriction}.

Finally, we check the condition (c). Let  $x\in X$ and $m\in\mathbb N$. If $x\in A_m$, then $h_{n+1,m}(x)=f_{n+1,m}(x)$, and if  $x\not\in A_m$, then $h_{n+1,m}(x)=h_{n,m}(x)$. In both cases 
$$
d(h_{n+1,m}(x),h_{n,m}(x))<\frac{1}{2^{n}}.
$$

Therefore, we have constructed sequences of mappings  $(h_{n,m})_{m=1}^\infty$ which satisfy (a)--(c) for every  $n\in\mathbb N$.

We prove that $h_{n,n}\mathop{\longrightarrow}\limits_{n\to\infty}f$ pointwise on $X$. Fix $x\in X$ and $\varepsilon>0$. Choose $n_0\in\mathbb N$ such that
$$
d(f(x),f_{n_0}(x))<\frac{1}{2^{n_0}}<\frac{\varepsilon}{2}.
$$
There exists a number $n_1\ge n_0$ such that 
$$
h_{n_0,n}(x)=f_{n_0}(x)
$$
for all $n\ge n_1$. Hence, for all $n\ge n_1$ we have
\begin{gather*}
  d(f(x),h_{n,n}(x))\le d(f(x),f_{n_0}(x))+\sum\limits_{k=n_0}^{n-1} d(h_{k,n}(x),h_{k+1,n}(x))<\frac{1}{2^{n_0}}+\frac{1}{2^{n_0}}<\varepsilon.
\end{gather*}
By the inductive assumption  $h_{n,n}\in \Lambda_{<\alpha}(X,Y)$ for every $n\in\mathbb N$. Hence, $f\in \Lambda_{\alpha}(X,Y)$. Consequently, the proposition (i) is proved for all $\alpha\le\omega_0$ and the proposition (ii) is proved for all limit ordinals $\alpha$.

It $\alpha=\beta+m$, where $\beta$ is a limit ordinal and $m\in\mathbb N$, and $f\in \Sigma_{\alpha+1}^f(X,Y)$, then, by the same method as in proof of  (i), one can show that there exists a sequence of mappings $g_n\in \Sigma^f_{\beta+m}(X,Y)$ which is pointwise convergent to $f$ on $X$. By the inductive assumption, $g_n\in \Lambda_{\beta+m-1}(X,Y)$, and hence $f\in \Lambda_{\alpha}(X,Y)$.
\end{pf}

\begin{thm}\label{closed_pointwise_sigma}
Let $X$ be a topological space, $Y$ be a space with a $\sigma$-disjoint base. Then the class $\Sigma^f(X,Y)$ is closed under pointwise limits.
\end{thm}

\begin{pf}
Let $(f_n)_{n=1}^\infty$ be a sequence of mappings $f_n\in \Sigma^f(X,Y)$ which converges pointwise to a mapping $f:X\to Y$. We show that $f\in\Sigma^f(X,Y)$.

Consider a $\sigma$-disjoint base $\mathcal V=\bigcup\limits_{m=1}^\infty \mathcal V_m$ of $Y$ and a $\sigma$-sfd base $\mathcal B_n=\bigcup\limits_{m=1}^\infty\mathcal B_{n,m}$ for $f_n$, $n\in\mathbb N$. Denote $\mathcal B=\bigcup\limits_{n=1}^\infty \mathcal B_n$ and for every $B\in\mathcal B$ we put 
  $$
  \widetilde{\mathcal V}_B=(V\in\mathcal V:\exists n\,\,|\,\, f_n(B)\subseteq V).
  $$
Notice that the family $\widetilde{\mathcal V}_B$ is at most countable, since every family $\mathcal V_m$ is disjoint. Enumerate the family $\widetilde{\mathcal V}_B$ in a sequence $(V_{B,k}:k\in\mathbb N)$. For all  $n,m,k\in\mathbb N$ we put
  $$
  \mathcal B_{n,m,k}=(B\cap f^{-1}(V_{B,k}):B\in\mathcal B_{n,m}).
  $$
Clearly, $\mathcal B_{n,m,k}$ is an sfd-family. It remains to prove that the family 
  $$
  \widetilde{\mathcal B}=\bigcup\limits_{n,m,k}\mathcal B_{n,m,k}
  $$
is a base for $f$. Let $V\in\mathcal V$ and $x\in f^{-1}(V)$. Take a number  $n$ such that $f_n(x)\in V$. Since $\mathcal B_n$ is a base for  $f_n$, there are $m\in\mathbb N$ and $B\in\mathcal B_{n,m}$ such that  $x\in B\subseteq f_n^{-1}(V)$. Then  $V\in\widetilde{\mathcal V}_B$. Hence, $x\in B\cap f^{-1}(V_{B,k})\subseteq f^{-1}(V)$.
\end{pf}

\begin{thm}\label{thm:pointwise_limit}
  Let  $X$ be a topological space, $Y$  be a perfectly normal space,  \mbox{$0\le\alpha<\omega_1$} and let $(f_n)_{n=1}^\infty$ be a sequence of mappings  $f_n\in K_{\alpha}(X,Y)$ which converges pointwise to a mapping $f:X\to Y$.  Then $f\in K_{\alpha+1}(X,Y)$.
\end{thm}

\begin{pf}
Let $F$ be a closed subset of $Y$ and $(G_n)_{n=1}^\infty$ be a decreasing sequence of open subsets of  $Y$ such that 
  $$
  F=\bigcap\limits_{n=1}^\infty G_n=\bigcap\limits_{n=1}^\infty \overline{G}_n.
  $$
It follows from the equality $\lim\limits_{n\to\infty}f_n(x)=f(x)$ for every $x\in A$ that 
  $$
  f^{-1}(F)=\bigcap\limits_{n=1}^\infty\bigcup\limits_{k=n+1}^\infty f_{k}^{-1}(G_n).
  $$
Since $f_k\in K_\alpha(X,Y)$, for every  $n$ the set $f_k^{-1}(G_n)$ belongs to the $\alpha$'th functionally additive class, hence, $f^{-1}(F)$ is a set of the $(\alpha+1)$'th functionally multiplicative class in $X$.
\end{pf}

\begin{thm}\label{thm:Lambda_is_Sigma}
Let $X$ be a topological space, $Y$ be a perfectly normal space with a $\sigma$-disjoint base and \mbox{$0<\alpha<\omega_1$}. Then
\begin{enumerate}
  \item[(i)] $\Lambda_{\alpha}(X,Y)\subseteq {\Sigma}_{\alpha}^f(X,Y)$, if $\alpha<\omega_0$,

  \item[(ii)] $\Lambda_{\alpha}(X,Y)\subseteq {\Sigma}_{\alpha+1}^f(X,Y)$, if $\alpha\ge \omega_0$.
\end{enumerate}
\end{thm}

\begin{pf}  The theorem is obvious for $\alpha=1$. Assume it is true for all $1\le\beta<\alpha$  and prove it for $\beta=\alpha$.
Let  $f\in \Lambda_\alpha(X,Y)$.

We consider the case $\alpha=\beta+1$ and take a sequence $(f_n)_{n=1}^\infty$ of mappings $f_n\in\Lambda_{\beta}(X,Y)$ which is pointwise convergent to збігається $f$ on $X$. By the assumption $f_n\in\Sigma_{\beta}^f(X,Y)$ in case $\alpha$ is finite, or $f_n\in\Sigma_{\alpha}(X,Y)$ in case  $\alpha$ is infinite. Then, respectively,  $f_n\in K_{\beta}(X,Y)\cap\Sigma^f(X,Y)$ or $f_n\in K_{\alpha}(X,Y)\cap\Sigma^f(X,Y)$ by Theorem~\ref{cor:sigma_F}. Applying Theorems~\ref{closed_pointwise_sigma} and~\ref{thm:pointwise_limit} we obtain that $f\in K_\alpha(X,Y)\cap\Sigma^f(X,Y)=\Sigma_\alpha^f(X,Y)$ if  $\alpha<\omega_0$, or $f\in K_{\alpha+1}(X,Y)\cap\Sigma^f(X,Y)=\Sigma_{\alpha+1}^f(X,Y)$, if $\alpha\ge\omega_0$.

Now we suppose  that $\alpha={\rm sup\,\,}\alpha_n$ is a limit ordinal and let  $(f_n)_{n=1}^\infty$ be a sequence of mappings $f_n\in\Lambda_{\alpha_n}(X,Y)$ which converges pointwise to $f$ on $X$. By the assumption $f_n\in\Sigma_{\alpha_n+1}^f(X,Y)\subseteq\Sigma_{\alpha}^f(X,Y)$ for every $n$. Theorems~\ref{cor:sigma_F},~\ref{closed_pointwise_sigma} and~\ref{thm:pointwise_limit} imply that $f\in K_{\alpha+1}^f(X,Y)\cap\Sigma^f(X,Y)=\Sigma_{\alpha+1}^f(X,Y)$.
\end{pf}

Combining Theorems~\ref{finite_classes} and \ref{thm:Lambda_is_Sigma} we get the following result.
\begin{thm}\label{thm:main_Banach}
  Let $X$ be a topological space, $Y$ be a metric space and  \mbox{$0<\alpha<\omega_1$}. Then
\begin{enumerate}
  \item[(i)] $\Lambda_{\alpha}(X,Y)={\Sigma}_{\alpha}^f(X,Y)$, if $\alpha<\omega_0$,

  \item[(ii)] $\Lambda_{\alpha}(X,Y)={\Sigma}_{\alpha+1}^f(X,Y)$, if $\alpha\ge \omega_0$.
\end{enumerate}
\end{thm}

Finally, we show that the condition on $X$ to be perfect in Theorem~\ref{thm:Hansell} is essential.

\begin{exm}
  {\it There exists a normal space $X$ such that ${\rm H}_2(X,\mathbb R)\ne \Phi_2(X,\mathbb R)$.}
\end{exm}

\begin{pf} Let  $A$ be a $G_{\delta\sigma}$-set which is not an $F_{\sigma\delta}$-set in $\mathbb R$,
and let \mbox{$X=(\mathbb R,\tau)$} be the real line with a topology $\tau$ such that a basic neighborhood of $x\in \mathbb R\setminus A$ is the set $\{x\}$ and a basic neighborhood of  $x\in A$ is an interval $(x-\varepsilon,x+\varepsilon)$, $\varepsilon>0$. The normality of $X$ follows from~\cite[Example 5.1.22]{Eng-eng}.

Since every open set in $\mathbb R$ is open in $X$,  $A$ is $G_{\delta\sigma}$ in $X$. Moreover, $A$ is closed in $X$, therefore, $A$ is an ambiguous set of the second class in $X$. Hence, the characteristic function $f=\chi_A:X\to\mathbb R$ of $A$ belongs to the class ${\rm H}_2(X,\mathbb R)$.

Notice that the normality of $X$ implies the equality ${\rm H}_1(X,\mathbb R)={\rm K_1}(X,\mathbb R)$ (see~\cite[Proposition 1.8]{Vesely}). Then $\Phi_2(X,\mathbb R)=\Lambda_2(X,\mathbb R)$.

We prove that $f\not\in\Phi_2(X,\mathbb R)$. To obtain a contradiction, suppose that $f\in\Lambda_2(X,\mathbb R)$. By Theorem~\ref{thm:main_Banach} there exists a sequence of functions $f_n\in\Sigma_1^f(X,\mathbb R)={\rm K_1}(X,\mathbb R)$ which converges pointwise to $f$ on $X$. Notice that every function $f_n:X\to\mathbb R$ is of the first Baire class  (see for instance \cite{Vesely}). Hence, $f:X\to\mathbb R$ is the function of the second Baire class. Let $(f_{n,m})_{n,m=1}^\infty$ be a sequence of continuous functions on  $X$ such that $f(x)=\lim\limits_{n\to\infty}\lim\limits_{m\to\infty}f_{n,m}(x)$ for every  $x\in X$. The definition of $\tau$ implies that for all $n,m$ the set $A$ is contained in the set $C(f_{n,m})$ of all points of continuity of  $f_{n,m}$ on $\mathbb R$. It is well-known that $C(f_{n,m})$ is a $G_\delta$-subset of $\mathbb R$. We put $B=\bigcap\limits_{n,m=1}^\infty C(f_{n,m})$. Then  $B$ is a  $G_\delta$-set in $\mathbb R$ which contains $A$ and the restriction  $f|_B$ belongs to the second Baire class in the Euclidean topology. Since $A=(f|_B)^{-1}((0,1])$, $A$ is an $F_{\sigma\delta}$-subset of $B$ and, consequently, is an $F_{\sigma\delta}$-subset of  $\mathbb R$, a contradiction.
\end{pf}

\end{document}